\documentstyle[11pt]{article}
\newcommand{\ba}{\begin{array}}
\newcommand{\ea}{\end{array}}
\newcommand{\be}{\begin{em}}
\newcommand{\ee}{\end{em}}
\newcommand{\bd}{\begin{document}}
\newcommand{\ed}{\end{document}}
\newcommand{\bn}{\begin{eqnarray}}
\newcommand{\en}{\end{eqnarray}}
\newcommand{\bq}{\begin{equation}}
\newcommand{\eq}{\end{equation}}
\newcommand{\al}{\alpha}
\newcommand{\ga}{\gamma}
\newcommand{\Ga}{\Gamma}
\newcommand{\de}{\delta}
\newcommand{\De}{\Delta}
\newcommand{\ve}{\varepsilon}
\newcommand{\la}{\lambda}
\newcommand{\La}{\Lambda}
\newcommand{\si}{\sigma}
\newcommand{\Si}{\Sigma}
\newcommand{\om}{\omega}
\newcommand{\Om}{\Omega}
\newcommand{\vp}{\varphi}
\newcommand{\ol}{\overline}
\newcommand{\ul}{\underline}
\newcommand{\ti}{\tilde}
\newcommand{\bs}{\backslash}
\newcommand{\pa}{\parallel}
\newcommand{\lf}{\left\{}
\newcommand{\ri}{\right\}}
\newcommand{\r}{\right}
\newcommand{\rb}{\rbrack}
\newcommand{\lb}{\lbrack}
\newcommand{\p}{\partial}
\newcommand{\R}{{\rm Re}}
\newcommand{\I}{{\rm Im}}
\newcommand{\n}{\eqno}
\newcommand{\vv}{\vspace{2.0mm}}
\newcommand{\vs}{\vspace{1mm}}
\newcommand{\vsp}{\vspace{2.0mm}}
\newcommand{\vsa}{\vspace{2.2mm}}
\newcommand{\vsc}{\vspace{1.8mm}}
\newcommand{\vse}{\vspace{1.4mm}}
\newcommand{\h}{\hfill}
\newcommand{\On}{{\rm on}}
\newcommand{\In}{{\rm in}}
\newcommand{\f}{{\rm for}}
\newcommand{\If}{{\rm if}}
\newcommand{\Or}{{\rm or}}
\newcommand{\bvp}{{\rm boundary\;value\;problem}}
\newcommand{\bc}{{\rm boundary\;condition}}
\newcommand{\ece}{{\rm elliptic\;complex\;equation}}
\newcommand{\ce}{{\rm complex\;equation}}
\newcommand{\af}{{\rm analytic\;function}}
\newcommand{\CI}{C\!\!\!\!I}
\newcommand{\IR}{I\!\!R}
\newcommand{\di}{\displaystyle}
\begin{document}
\title{\bf The equivalent theorem of a new generalized Bernstein-B$\acute{e}$zier operators
}
\author{\small Qiu-Lan Qi$^{1,2}$\thanks{Corresponding authors: Qiu-Lan Qi and Dan-Dan Guo   }, ~ Dan-Dan Guo$^{1\ast}$, ~ Ge Yang$^{1,2}$  \\
{\small\sl $^{1}$School of Mathematical Sciences,
 Hebei Normal University,} \\
{\small\sl Shijiazhuang 050024, China}\\
{\small\sl $^{2}$Hebei Key Laboratory of Computational Mathematics and Applications,} \\
{\small\sl Shijiazhuang 050024, China }\\
{\small\sl E-mail: qiqiulan@163.com, 2728561580@qq.com, yanggeshida@163.com}}\date{}
\maketitle{\bf Abstract}. In this paper, a new generalized Bernstein-B$\acute{e}$zier type operators is constructed. The estimates of the moments of these operators are investigated. The rate of convergence in terms of modulus of continuity is given. Then, the equivalent theorem of these operators is studied.
{\bf 1.\,\,\,\,Introduction}\\
The famous Bernstein operators are defined by$^{[1]}$
$$B_n(f;x)=\sum\limits^{n}_{k=0}f(\frac{k}{n})p_{n,k}(x),$$
where  $f(x)\in C[0,1] , p_{n,k}(x)=\left(^{n}_{k}\right)x^{k}(1-x)^{n-k}, k=1,2....n.$\\
The Bernstein-B$\acute{e}$zier polynomials for any $\alpha >0,$ and a function $f$ defined on $[0,1]$ as follows$^{[2]}$:
$$B_{n,\alpha}(f,x)=\sum^{n}_{k=0}f(\frac{k}{n})\left[J^{\alpha}_{n,k}(x)-J^{\alpha}_{n,k+1}(x)\right] ,\eqno{(1)}$$
where $J_{n,k}(x)=\sum^{n}_{i=k}p_{n,i}(x), k=0,1,....n$, $J_{n,k}(x)$ is the B$\acute{e}$zier basis function of degree $n$, and
 $1>J_{n,0}(x)>J_{n,1}(x)>\cdot\cdot\cdot>J_{n,n}(x)=x^n$.\\
In Computer Aided Geometric Design, B$\acute{e}$zier
basis functions are very useful and their analytical properties have been studied by many authors$^{[2-8]}$. In 2013, Ren introduced Bernstein operators as follows$^{[8]}$: for $f(x)\in C[0,1], \beta\in[0,1]$,
$$E_{n,\beta}(f;x)=f(0)p_{n,0}(x)+\sum^{n-1}_{k=1}p_{n,k}(x)F_{n,k}^{(\beta)}(f)+f(1)p_{n,n}(x),$$
where
$$F_{n,k}^{(\beta)}(f)=\frac{1}{B(nk,n(n-k))}\int^{1}_{0}t^{nk-1}(1-t)^{n(n-k)-1}f(\beta t+(1-\beta)\frac{k}{n})dt, $$
 $B(\cdot ,\cdot)$ is the Beta function.\\
In 2017, Ren$^{[9]}$ studied generalized Bernstein operators $E_{n,\beta}^{\alpha}(f;x)$, and obtained the direct theorem for these operators. In order to investigate the inverse theorem, using $p_{n-1,k}(x),$ we redefine the generalized Bernstein-B$\acute{e}$zier operators as follows: 
For $\beta\in[0,1], \alpha\geq1,$
$$ L_{n,\beta}^{(\alpha)}(f;x)=f(0)Q_{n,0}^{(\alpha)}(x)+\sum^{n-1}_{k=1}Q_{n,k}^{(\alpha)}(x)F_{n-1,k}^{(\beta)}(f)+f(1)Q_{n,n}^{(\alpha)}(x) ,$$
where $Q_{n,k}^{(\alpha)}(x)=J^{\alpha}_{n,k}(x)-J^{\alpha}_{n,k+1}(x),$ $J_{n,k}(x)$ and $F_{n-1,k}^{(\beta)}(f)$ are defined as above.\\
The operators $L_{n,\beta}^{(\alpha)}(f;x)$ are bounded and positive on C$[0,1]$. When $\alpha=1$, $L_{n,\beta}^{(\alpha)}(f;x)$ become the operators $L_{n,\beta}(f;x)$. When $\beta=0$, $L_{n,\beta}^{(\alpha)}(f;x)$ become the Bernstein-B$\acute{e}$zier operators $B_{n,\alpha}(f,x)$ (see (1)) .\\
{\bf Remark 1}  Throughout this paper, $M$ is a positive constant independent of $n$ and $x$, the value of $M$ may be different in different places.\\
{\bf Remark 2}  In this paper, for $f\in C[0,1]$:=$\{f:f$ is continous on $[0,1]\}$ , the norm of $f(x)$ is defined as
 $ \|f\|=\max\{|f(x)|: x\in [0,1]\}.$\\
{\bf Remark 3}$^{[1]}$   For $B_n(t^{j};x), j=0,1,2,$ one has
\begin{eqnarray*}
&(1)&B_n(1;x)=1;\\
&(2)&B_n(t;x)=x;\\
&(3)&B_n(t^{2};x)=x^{2}+\frac{x(1-x)}{n}.\\
\end{eqnarray*}
{\bf 2.\,\,\,\,Estimates of moments }\\
In this section, using similar calculation method in reference [2, 7, 9], we can obtain the estimate of moment, here we omit the calculation details.\\
{\bf Lemma 1}$^{[9\  Lemma 2]}$  For $F_{n-1,k}^{(\beta)}(t^{i})$, $i=0, 1, 2,$  we have
\begin{eqnarray*}
&(1)&F_{n-1,k}^{(\beta)}(1)=1;\\
&(2)&F_{n-1,k}^{(\beta)}(t)=\frac{k}{n-1};\\
&(3)&F_{n-1,k}^{(\beta)}(t^{2})=\frac{\beta^{2}}{(n-1)^{2}+1}\cdot \frac{k}{n-1}+\left(1-\frac{\beta^{2}}{(n-1)^{2}+1}\right)\cdot\frac{k^{2}}{(n-1)^{2}}.
\end{eqnarray*}
{\bf Lemma 2} For $L_{n,\beta}(t^{i};x)$, $i=0,1,2,$ we have
\begin{eqnarray*}
&(1)&L_{n,\beta}(1;x)=1;\\
&(2)&L_{n,\beta}(t;x)=x+\frac{1}{n-1}(x-x^{n});\\
&(3)&L_{n,\beta}(t^{2};x)=\frac{n^{2}}{(n-1)^{2}}\cdot x^{2}+[\frac{n}{(n-1)^{2}}+\frac{\beta^{2}}{(n-1)^{2}+1}\cdot \frac{n}{n-1}]\cdot \varphi^{2}(x)\\
    &&-\frac{\beta^{2}}{(n-1)^{2}+1}\cdot \frac{n}{(n-1)^{2}} \cdot x+[\frac{\beta^{2}}{(n-1)^{2}+1}\cdot \frac{n}{(n-1)^{2}}-\frac{2n-1}{(n-1)^{2}}]\cdot x^{n}.
\end{eqnarray*}
{\bf Lemma 3} Let $\alpha\geq1$, we have
\begin{eqnarray*}
&(1)&\frac{1}{n-1}\sum^{n-1}_{k=1}J_{n,k}(x)=\frac{n}{n-1}x-\frac{1}{n-1}x^{n};\\
&(2)&\frac{1}{(n-1)^{2}}\sum^{n-1}_{k=1}kJ_{n,k}(x)=\frac{n}{2(n-1)}x^{2}.
\end{eqnarray*}
{\bf Lemma 4}$^{[2,\  9\  Lemma 1]}$  Let $\alpha\geq1$, we have
\begin{eqnarray*}
&(1)&\lim_{n\rightarrow\infty}\frac{1}{n-1}\sum^{n-1}_{k=1}J^{\alpha}_{n,k}(x)=x \hspace{0.3cm} uniformly \hspace{0.2cm} on \hspace{0.1cm} [0,1];\\
&(2)&\lim_{n\rightarrow\infty}\frac{1}{(n-1)^{2}}\sum^{n-1}_{k=1}kJ^{\alpha}_{n,k}(x)=\frac{x^{2}}{2}   \hspace{0.3cm} uniformly  \hspace{0.2cm} on  \hspace{0.1cm} [0,1].
\end{eqnarray*}
{\bf Lemma 5}$^{[7,\  9\  Lemma 3]}$  For $x\in[0,1]$, $\alpha\geq1$, $k=0,1,...n-1,$ we have
$$0\leq Q_{n,k}^{(\alpha)}(x)\leq \alpha P_{n,k}(x).$$
{\bf Lemma 6} Let $\alpha\geq1$, we have
\begin{eqnarray*}
&(1)&L_{n,\beta}^{(\alpha)}(1;x)=1;\\
&(2)&\lim_{n\rightarrow\infty}L_{n,\beta}^{(\alpha)}(t;x)=x \hspace{0.3cm} uniformly \hspace{0.2cm} on \hspace{0.1cm} [0,1];\\
&(3)&\lim_{n\rightarrow\infty}L_{n,\beta}^{(\alpha)}(t^{2};x)=x^{2} \hspace{0.3cm} uniformly \hspace{0.2cm} on \hspace{0.1cm} [0,1];
\end{eqnarray*}
{\bf Lemma 7} For $\alpha\geq1$, $\beta\in [0,1], \varphi^{2}(x)=x(1-x)$, we have
\begin{eqnarray*}
&(1)&L_{n,\beta}^{(\alpha)}((t-x)^{2};x)\leq\frac{\alpha}{4}(14+\beta^{2})\cdot \frac{1}{n-1},  \hspace{0.4cm}for \hspace{0.2cm} x\in[0,1];\\
&(2)&L_{n,\beta}^{(\alpha)}((t-x)^{2};x)\leq \frac{3\varphi^{2}(x)}{n-1},  \hspace{0.4cm} for\hspace{0.2cm} x\in E_{n}=\left[\frac1n, 1-\frac1n\right].
\end{eqnarray*}
{\bf Lemma 8} Let $f(x)\in C[0,1]$, $\varphi(x)=\sqrt{x(1-x)}$, $0\leq\lambda\leq1$, we have
$$|\varphi^{\lambda}(x)\cdot(L_{n,\beta}^{(\alpha)}(f;x))'|\leq 15\alpha\varphi^{\lambda-1}(x)\sqrt{n}\|f\| .$$
{\bf Lemma 9} Let $f\in W_{\lambda}$, $\varphi(x)=\sqrt{x(1-x)}$, we have
$$|\varphi^{\lambda}(x)\cdot(L_{n,\beta}^{(\alpha)}(f;x))'|\leq 104\alpha\|\varphi^{\lambda}f'\| .$$
{\bf 3.\,\,\,\,Direct Theorem}\\
For $f(x)\in C[0,1]$, $\varphi(x)=\sqrt{x(1-x)}$, $\lambda\in[0,1]$, $x\in [0,1]$, let
$$\omega_{\varphi^{\lambda}}(f;t)=\sup_{0< h\leq t} \sup_{x\pm\frac{h\varphi^{\lambda}(x)}{2}\in[0,1]}\left|f(x+\frac{h\varphi^{\lambda}(x)}{2})-f(x-\frac{h\varphi^{\lambda}(x)}{2})\right|,$$
be the Ditzian-Totik modulus, and let
$$K_{\varphi^{\lambda}}(f;t)=\inf_{g\in W_{\lambda}} \{||f-g||+t||\varphi^{\lambda} g^{'}||\},$$
be the corresponding K-functional, here $W_{\lambda}=\{f|f\in A.C._{loc.}[0,1],\|\varphi^{\lambda} f^{'}\|<\infty, \|f^{'}\|<\infty \}.$
It is well know that$^{[1]}$
$$K_{\varphi^{\lambda}}(f;t)\sim \omega_{\varphi^{\lambda}}(f;t).$$
{\bf Theorem 1 (Korovkin Theorem)} Let $f(x)\in C[0,1]$, $\alpha\geq1$, then Bernstein-B$\acute{e}$zier operators $L_{n,\beta}^{(\alpha)}(f;x)$ converge to $f(x)$ on $[0,1]$.\\
{\bf Theorem 2} For $f(x)\in C[0,1]$, $\alpha\geq1$, $\varphi(x)=\sqrt{x(1-x)}$ , $0\leq\beta, \lambda\leq1, $
then we have
$$|L_{n,\beta}^{(\alpha)}(f;x)-f(x)|\leq M\omega_{\varphi^{\lambda}}(f;\frac{\varphi^{1-\lambda}(x)}{\sqrt{n}}). $$
{\sl Proof} \hspace{0.3cm} Let $g\in W_{\lambda} $, then
 $|L_{n,\beta}^{(\alpha)}(f;x)-f(x)|\leq 2\|f-g\|+|L_{n,\beta}^{(\alpha)}(g;x)-g(x)|.$\\
Since $g(t)=\int^{t}_{x}g'(u)du+g(x),$  $L_{n,\beta}^{(\alpha)}(1;x)=1,$
we know
$$
|L_{n,\beta}^{(\alpha)}(g;x)-g(x)|\leq\|\varphi^{\lambda}g'\|L_{n,\beta}^{(\alpha)}\left(|\int^{t}_{x}\varphi^{-\lambda}(u)du|;x\right).
$$
By the H\"{o}lder inequality, one has
$$
\left|\int^{t}_{x}\varphi^{-\lambda}(u)du\right|\leq 4\varphi^{-\lambda}(x)\cdot|t-x|,\eqno{(2)}
$$
 we get
 $|L_{n,\beta}^{(\alpha)}(g;x)-g(x)|\leq \|\varphi^{\lambda}g'\|\cdot\varphi^{-\lambda}(x)\cdot L_{n,\beta}^{(\alpha)}(|t-x|;x).$\\
By  the definition of $L_{n,\beta}^{(\alpha)}(f;x)$ and Lemma 7 (2), we have
$$|L_{n,\beta}^{(\alpha)}(g;x)-g(x)|_{E_{n}}\leq M\|\varphi^{\lambda}g'\|\cdot\frac{\varphi^{1-\lambda}(x)}{\sqrt{n}}.$$
Combing Theorem $8.4.8^{[1]}$, we have
$$|L_{n,\beta}^{(\alpha)}(f;x)-f(x)|\leq M\left(\|f-g\|+\|\varphi^{\lambda}g'\|\cdot\frac{\varphi^{1-\lambda}(x)}{\sqrt{n}}\right).
$$
Taking infimum on the right hand side over all $g\in W_{\lambda},$ we get the desired result.\\
{\bf Theorem 3}  For $f(x)\in C^{1}[0,1]$:=$\{f:f'(x)$ is continous on $[0,1]\}$ , and $\alpha\geq1$, then
$$|L_{n,\beta}^{(\alpha)}(f;x)-f(x)|
\leq\sqrt{\frac{(14+\beta^{2})\alpha}{4n}}\cdot\left\{\|f'\|+\omega(f';\frac{1}{\sqrt{n}})\cdot
(1+\sqrt{\frac{(14+\beta^{2})\alpha}{4}})\right\}.
$$
{\sl Proof} \hspace{0.3cm}For any $t, x\in[0,1]$, $\delta>0 $, by the Taylor's expansion, we get
$$
|f(t)-f(x)-f'(x)(t-x)|\leq\left|\int^{t}_{x}|f'(u)-f'(x)|du\right|
\leq\omega(f';\delta)(|t-x|+\delta^{-1}(t-x)^{2}),
$$
hence, by the Cauchy-Schwarz inequality, we have
\begin{eqnarray*}
&&|L_{n,\beta}^{(\alpha)}(f(t)-f(x)-f'(x)(t-x);x)|\\
&\leq& \omega(f';\delta)(L_{n,\beta}^{(\alpha)}(|t-x|;x)+\delta^{-1}L_{n,\beta}^{(\alpha)}((t-x)^{2});x)\\
&\leq&\omega(f';\delta)\cdot\left[\sqrt{L_{n,\beta}^{(\alpha)}(1;x)}+\delta^{-1}\sqrt{L_{n,\beta}^{(\alpha)}((t-x)^{2};x)}\right]\cdot\sqrt{L_{n,\beta}^{(\alpha)}((t-x)^{2};x)}.
\end{eqnarray*}
Thus,
\begin{eqnarray*}
&&|L_{n,\beta}^{(\alpha)}(f;x)-f(x)|\\
&\leq&\|f'\|\cdot L_{n,\beta}^{(\alpha)}(|t-x|;x)+
\omega(f';\delta)\cdot\left[1+\delta^{-1}\sqrt{L_{n,\beta}^{(\alpha)}((t-x)^{2};x)}
\right]\cdot\sqrt{L_{n,\beta}^{(\alpha)}((t-x)^{2};x)}.
\end{eqnarray*}
Taking $\delta=\frac{1}{\sqrt{n}}$, by Lemma 7 (1), we can obtain Theorem 3.\\
{\bf 4.\,\,\,\,Equivalent  Theorem}\\
{\bf Theorem 4} Let $f(x)\in C[0,1]$, $\varphi(x)=\sqrt{x(1-x)}$, $  0\leq\beta, \gamma, \lambda\leq1, \alpha\geq1 $, if
$$|L_{n,\beta}^{(\alpha)}(f;x)-f(x)|=O(n^{-\frac{\gamma}{2}}),$$
one has $\omega_{\varphi^{\lambda}}(f;t)=O(t^{\gamma}).$\\
{\sl Proof} \hspace{0,3cm}By the definition of the $K$-functional, 
\begin{eqnarray*}
K_{\varphi^{\lambda}}(f;t)&\leq&\|f-L_{n,\beta}^{(\alpha)}(f;x)\|+t\|\varphi^{\lambda}(x) (L_{n,\beta}^{(\alpha)}(f;x))'\|\\
&\leq& Mn^{-\frac{\gamma}{2}}+t(\|\varphi^{\lambda}(x) (L_{n,\beta}^{(\alpha)}((f-g);x))'\|+\|\varphi^{\lambda}(x)(L_{n,\beta}^{(\alpha)}(g;x))'\|\\
&\leq& Mn^{-\frac{\gamma}{2}}+t\sqrt{n}(\|f-g\|+\frac{1}{\sqrt{n}}\|\varphi^{\lambda}g'\|)\\
&\leq& M(n^{-\frac{\gamma}{2}}+t\sqrt{n}\cdot K_{\varphi^{\lambda}}(f;n^{-\frac{1}{2}})).
\end{eqnarray*}
Applying the Berens-Lorentz Lemma, and the relation $\omega_{\varphi^{\lambda}}(f;t)\sim K_{\varphi^{\lambda}}(f;t),$ we obtain the desired Theorem 4.\\
From Theorem 2 and Theorem 4, we can obtain the equivalent theorem.\\
{\bf Theorem 5} Let $f(x)\in C[0,1]$, $\varphi(x)=\sqrt{x(1-x)}$, $0\leq\beta, \gamma, \lambda\leq1, \alpha\geq1 $, we have
$$|L_{n,\beta}^{(\alpha)}(f;x)-f(x)|=O(n^{-\frac{\gamma}{2}}) \Leftrightarrow \omega_{\varphi^{\lambda}}(f;t)=O(t^{\gamma}). $$
{\bf 5.\,\,\,\,Proof of Lemma 8 and Lemma 9}\\
{\sl Proof of Lemma 8:}  We write
$$\left(L_{n,\beta}^{(\alpha)}(f;x)\right)'=f(0)(Q_{n,0}^{(\alpha)}(x))'+\left(\sum^{n-1}_{k=1}Q_{n,k}^{(\alpha)}(x)F_{n-1,k}^{(\beta)}(f)\right)'+f(1)(Q_{n,n}^{(\alpha)}(x))'
=I_{1}+I_{2}+I_{3},\eqno{(3)}$$
and will estimate $I_{1}$ , $I_{2}$ and $I_{3}$, respectively. \\
Firstly, noting that $J_{n,0}'(x)=0,$ we have
\begin{eqnarray*}
|I_{1}|&=&|\alpha f(0)J_{n,0}^{\alpha-1}(x)\cdot J_{n,0}'(x)-\alpha f(0)J_{n,1}^{\alpha-1}(x)\cdot J_{n,1}'(x)|\\
&=&|n\alpha f(0)\cdot [1-(1-x)^{n}]^{\alpha-1}\cdot (1-x)^{n-1}|\\
&\leq& n\alpha\|f\|=\alpha\sqrt{n}\|f\|\cdot\sqrt{n}.
\end{eqnarray*}
For $x\in [0,\frac{1}{n}),$ one has
$$\varphi^{\lambda}(x)\cdot |I_{1}|\leq \alpha\varphi^{\lambda-1}\sqrt{n}\|f\|\cdot\sqrt{n}\cdot\sqrt{x(1-x)}\leq\alpha\varphi^{\lambda-1}\sqrt{n}\|f\|
\cdot\sqrt{n}\cdot\sqrt{x}\leq\alpha\varphi^{\lambda-1}(x)\sqrt{n}\|f\|.\eqno{(4)}$$
For $x\in (1-\frac{1}{n},1],$ one has
$$\varphi^{\lambda}(x)\cdot |I_{1}|\leq \alpha\varphi^{\lambda-1}\sqrt{n}\|f\|\cdot\sqrt{n}\cdot\sqrt{x(1-x)}\leq\alpha\varphi^{\lambda-1}\sqrt{n}\|f\|
\cdot\sqrt{n}\cdot\sqrt{1-x}\leq\alpha\varphi^{\lambda-1}(x)\sqrt{n}\|f\|.\eqno{(5)}$$
For $x\in E_{n}=[\frac{1}{n},1-\frac{1}{n}],$ one has
$$|I_{1}|=|-\alpha f(0)J_{n,1}^{\alpha-1}(x)J_{n,1}'(x)|\leq|-\alpha f(0) J_{n,1}'(x)|\leq|-\alpha f(0)\cdot[1-p_{n,0}(x)]'|=|\alpha f(0)\cdot p_{n,0}'(x)|
.$$
Combining that $$p_{n,k}'(x)=\frac{n}{\varphi^{2}(x)}\left(\frac{k}{n}-x\right)\cdot p_{n,k}(x),\hspace {0.5cm}B_{n}((t-x)^{2};x)=\frac{\varphi^{2}(x)}{n},$$
we get
\begin{eqnarray*}
\varphi^{\lambda}(x)\cdot |I_{1}|&\leq&\alpha\|f\|\cdot|p_{n,0}'(x)|\cdot\varphi^{\lambda-1}(x)\\
&\leq&\alpha\|f\|\cdot\sum^{n}_{k=0}\left|p_{n,k}'(x)\right|\cdot\varphi^{\lambda-1}(x)\\
&\leq&\alpha\|f\|\cdot\frac{n}{\varphi^{2}(x)}\sum^{n}_{k=0}\left|\frac{k}{n}-x\right|p_{n,k}(x)\cdot\varphi^{\lambda-1}(x)\\
&\leq&\alpha\|f\|\cdot\frac{n}{\varphi^{2}(x)}\cdot\sqrt{\sum^{n}_{k=0}\left(\frac{k}{n}-x\right)^{2}p_{n,k}(x)}\cdot\left(\sum^{n}_{k=0}p_{n,k}(x)\right)^{\frac{1}{2}}\cdot\varphi^{\lambda-1}(x)\\
&\leq&\alpha\varphi^{\lambda-1}(x)\sqrt{n}\|f\|.\hspace{9.7cm}(6)
\end{eqnarray*}
From (4)-(6), one has $\varphi^{\lambda}(x)\cdot |I_{1}|\leq\alpha\varphi^{\lambda-1}(x)\sqrt{n}\|f\|.$\\
Secondly, noting that $J_{n,n+1}'(x)=0,$ we have
\begin{eqnarray*}
|I_{3}|&=&|\alpha f(1)J_{n,n}^{\alpha-1}(x)\cdot J_{n,n}'(x)-\alpha f(1)J_{n,n+1}^{\alpha-1}(x)\cdot J_{n,n+1}'(x)|\\
&=&|n\alpha f(1)\cdot (x^{n})^{\alpha-1}\cdot x^{n-1}|\\
&\leq& n\alpha\|f\|=\alpha\sqrt{n}\|f\|\cdot\sqrt{n}.
\end{eqnarray*}
For $x\in [0,\frac{1}{n})\bigcup(1-\frac{1}{n},1],$ by the similar method used in (4) and (5), we can get
$$\varphi^{\lambda}(x)\cdot |I_{3}|\leq\alpha\varphi^{\lambda-1}(x)\sqrt{n}\|f\|.\eqno{(7)}$$
For $x\in E_{n}=[\frac{1}{n},1-\frac{1}{n}],$ combining $|p_{n,n}'(x)|\leq \sum^{n}_{k=0}|p_{n,k}'(x)|$, by the similar method used in (6), we can obtain (7).\\
Finally, considering $I_{2}$, recalling that $|F_{n-1,k}^{(\beta)}(f)| \leq \|f\|,$ we write
\begin{eqnarray*}
I_{2}&=&\alpha \sum^{n-1}_{k=1}F_{n-1,k}^{(\beta)}(f)J_{n,k}^{\alpha-1}(x)\cdot J_{n,k}'(x)-
\alpha \sum^{n-1}_{k=1}F_{n-1,k}^{(\beta)}(f)J_{n,k+1}^{\alpha-1}(x)\cdot J_{n,k+1}'(x)\\
&=&\alpha \sum^{n-1}_{k=1}F_{n-1,k}^{(\beta)}(f)J_{n,k}^{\alpha-1}(x)\cdot J_{n,k+1}'(x)+
\alpha \sum^{n-1}_{k=1}F_{n-1,k}^{(\beta)}(f)J_{n,k}^{\alpha-1}(x)\cdot p_{n,k}'(x)\\
&&-\alpha \sum^{n-1}_{k=1}F_{n-1,k}^{(\beta)}(f)J_{n,k+1}^{\alpha-1}(x)\cdot J_{n,k+1}'(x)\\
&=&\alpha \sum^{n-1}_{k=1}F_{n-1,k}^{(\beta)}(f)\left\{[J_{n,k}^{\alpha-1}(x)-J_{n,k+1}^{\alpha-1}(x)]J_{n,k+1}'(x)+J_{n,k}^{\alpha-1}(x)\cdot p_{n,k}'(x)\right\},
\end{eqnarray*}
so
$$
|I_{2}|\leq\alpha\|f\|\left(\sum^{n-1}_{k=1}[J_{n,k}^{\alpha-1}(x)
-J_{n,k+1}^{\alpha-1}(x)]J_{n,k+1}'(x)+\sum^{n-1}_{k=1}J_{n,k}^{\alpha-1}(x)\cdot |p_{n,k}'(x)|\right)
=\alpha\|f\|(Q_{1}+Q_{2}).
$$
Noting that $J_{n,1}'(x)>0, J_{n,0}'(x)=0, J_{n,n+1}(x)=0$, one can get
\begin{eqnarray*}
Q_{1}&=&\sum^{n-1}_{k=1}J_{n,k}^{\alpha-1}(x)J_{n,k+1}'(x)-\sum^{n-1}_{k=1}J_{n,k+1}^{\alpha-1}(x)J_{n,k+1}'(x)\\
&=&\sum^{n-1}_{k=1}J_{n,k}^{\alpha-1}(x)J_{n,k}'(x)-\sum^{n-1}_{k=1}J_{n,k}^{\alpha-1}(x)p_{n,k}'(x)-\sum^{n-1}_{k=1}J_{n,k+1}^{\alpha-1}(x)J_{n,k+1}'(x)\\
&=&-\sum^{n-1}_{k=1}J_{n,k}^{\alpha-1}(x)\cdot p_{n,k}'(x)-J_{n,n}^{\alpha-1}(x)p_{n,n}'(x)+J_{n,1}^{\alpha-1}(x)J_{n,1}'(x),\\
&\leq&\sum^{n-1}_{k=1}J_{n,k}^{\alpha-1}(x)\cdot |p_{n,k}'(x)|+J_{n,n}^{\alpha-1}(x)|p_{n,n}'(x)|+J_{n,1}^{\alpha-1}(x)J_{n,1}'(x)\\
&\leq&Q_{2}+p_{n,n}'(x)+[1-p_{n,0}(x)]'\\
&\leq&Q_{2}+p_{n,n}'(x)+|p_{n,0}'(x)|.
\end{eqnarray*}
Since
$$\alpha\|f\|\cdot p_{n,n}'(x)=\alpha\|f\|\cdot n\cdot x^{n-1}\leq\alpha\|f\|\sqrt{n}\cdot\sqrt{n},$$
$$\alpha\|f\|\cdot |p_{n,0}'(x)|=\alpha\|f\|\cdot n\cdot (1-x)^{n-1}\leq\alpha\|f\|\sqrt{n}\cdot\sqrt{n},$$
simimar to the estimate of $I_{1} $, we can obtain
$$\varphi^{\lambda}(x)\cdot\alpha\|f\|\cdot p_{n,n}'(x)\leq\alpha\varphi^{\lambda-1}(x)\sqrt{n}\|f\|,        $$
$$\varphi^{\lambda}(x)\cdot\alpha\|f\|\cdot |p_{n,0}'(x)|\leq\alpha\varphi^{\lambda-1}(x)\sqrt{n}\|f\|.        $$
Next, we think about $Q_{2}$.\\
In the case of $x=0$, one has
$$
p_{n,k}'(0)=\left\{\begin{array}{ll}
n,& for  \hspace{0.3cm} k=1;\\
0,& for  \hspace{0.3cm} 2\leq k\leq n-1.
\end{array}\right.
$$
In the case of $x=1$, one has
$$
p_{n,k}'(1)=\left\{\begin{array}{ll}
0,& for  \hspace{0.3cm} 1\leq k\leq n-2;\\
-n,& for  \hspace{0.3cm} k=n-1.
\end{array}\right.
$$
In the case of $0< x< 1, $ one has
$$p_{n,k}'(x)=\frac{n}{\varphi^{2}(x)}\left(\frac{k}{n}-x\right)\cdot p_{n,k}(x).$$
So
\begin{eqnarray*}
\varphi^{\lambda}(x)Q_{2}
&=&\varphi^{\lambda}(x)\left[J_{n,1}^{\alpha-1}(x)|p_{n,1}'(x)|+\sum^{n-2}_{k=2}J_{n,k}^{\alpha-1}(x)\cdot |p_{n,k}'(x)|+J_{n,n-1}^{\alpha-1}(x)|p_{n,n-1}'(x)|\right]\\
&\leq&\varphi^{\lambda}(x)|p_{n,1}'(x)|+\varphi^{\lambda}(x)\sum^{n-2}_{k=2}J_{n,k}^{\alpha-1}(x)\cdot |p_{n,k}'(x)|+\varphi^{\lambda}(x)|p_{n,n-1}'(x)|\\
&\leq&\varphi^{\lambda}(x)\left|n+\frac{\sqrt{n}}{\varphi(x)}+0\right|+\sum^{n-2}_{k=2}J_{n,k}^{\alpha-1}(x)|p_{n,k}'(x)|\varphi^{\lambda}(x)+\varphi^{\lambda}(x)\left|0-n+\frac{\sqrt{n}}{\varphi(x)}\right|\\
&=&\varphi^{\lambda}(x)n+\varphi^{\lambda-1}(x)\sqrt{n}+\varphi^{\lambda}(x) n+\varphi^{\lambda-1}(x)\sqrt{n}+\sum^{n-2}_{k=2}J_{n,k}^{\alpha-1}(x)|p_{n,k}'(x)|\varphi^{\lambda}(x)\\
&\leq& 4\varphi^{\lambda-1}(x)\sqrt{n}+\frac{n}{\varphi(x)}\sum^{n-2}_{k=2}\left|\frac{k}{n}-x\right|\cdot p_{n,k}(x)\cdot\varphi^{\lambda}(x) \\
&\leq& 4\varphi^{\lambda-1}(x)\sqrt{n}+\frac{n}{\varphi(x)}{\varphi^{\lambda}(x)}\cdot \left(B_{n}((t-x)^{2};x)\right)^{\frac{1}{2}}\cdot(p_{n,k}(x))^{\frac{1}{2}}\\
&\leq& 5\varphi^{\lambda-1}(x)\sqrt{n}.
\end{eqnarray*}
Then,
$$\varphi^{\lambda}(x)Q_{2}\leq 5\varphi^{\lambda-1}(x)\sqrt{n},\eqno{(8)}$$
and $\varphi^{\lambda}(x)Q_{1}\leq 5\varphi^{\lambda-1}(x)\sqrt{n}+2\sqrt{n}\cdot \varphi^{\lambda-1}(x)\cdot\frac{\sqrt{nx}}{e^{nx}}\leq 7\varphi^{\lambda-1}(x)\sqrt{n}.$\\
So
$$\varphi^{\lambda}(x)|I_{2}|\leq \alpha\|f\|\cdot \varphi^{\lambda}(x)\left(Q_{1}+Q_{2}\right)\leq 12\alpha\sqrt{n}\|f\|\cdot \varphi^{\lambda-1}(x).\eqno{(9)}$$
From (3)-(9), the desired result follows immediately.\\
{\sl  Proof of Lemma 9:}  Since $f(x)\left(L_{n,\beta}^{(\alpha)}(1;x)\right)'=0,$ we get
\begin{eqnarray*}
\left(L_{n,\beta}^{(\alpha)}(f;x)\right)'&=&f(0)\left(Q_{n,0}^{(\alpha)}(x)\right)'+\sum^{n-1}_{k=1}F_{n-1,k}^{(\beta)}(f)\left[J_{n,k}^{\alpha}(x)-J_{n,k+1}^{\alpha}(x)\right]'+f(1)\left(Q_{n,n}^{(\alpha)}(x)\right)'\\
&&-f(x)\left\{\left(Q_{n,0}^{(\alpha)}(x)\right)'+\sum^{n-1}_{k=1}F_{n-1,k}^{(\beta)}(1)\left[J_{n,k}^{\alpha}(x)-J_{n,k+1}^{\alpha}(x)\right]'+\left(Q_{n,n}^{(\alpha)}(x)\right)'\right\}\\
&=&[f(0)-f(x)]\left(Q_{n,0}^{(\alpha)}(x)\right)'+\sum^{n-1}_{k=1}\left[F_{n-1,k}^{(\beta)}(f)-f(x)\right]
\left[J_{n,k}^{\alpha}(x)-J_{n,k+1}^{\alpha}(x)\right]'\\
&&+[f(1)-f(x)]\left(Q_{n,n}^{(\alpha)}(x)\right)',
\end{eqnarray*}
we write
$$\left(L_{n,\beta}^{(\alpha)}(f;x)\right)'=H_{1}+H_{2}+H_{3},\eqno{}(10)$$
and will estimate $H_{1}$, $H_{2}$ and $H_{3}$ respectively. First, from (2), one can get
\begin{eqnarray*}
\varphi^{\lambda}(x)|H_{1}|&=&\varphi^{\lambda}(x)\cdot\left|\int^{0}_{x}\frac{\varphi^{\lambda}(u)f'(u)}{\varphi^{\lambda}(u)}du\right|\cdot
\left|\left(J_{n,0}^{\alpha}(x)\right)'-\left(J_{n,1}^{\alpha}(x)\right)'\right|\\
&\leq&\varphi^{\lambda}(x)\cdot\|\varphi^{\lambda}f'\|\cdot\left|\int^{0}_{x}\varphi^{-\lambda}(u)du\right|\cdot
\left|0-n\alpha (1-x)^{n-1}\right|\\
&\leq&4\varphi^{\lambda}(x)\cdot\|\varphi^{\lambda} f'\|\cdot \frac{|x|}{\varphi^{\lambda}(x)}\cdot n\alpha (1-x)^{n-1}\\
&\leq& 4\alpha\|\varphi^{\lambda} f'\|\cdot nx(1-x)^{n-1}.
\end{eqnarray*}
For $x\in[0,\frac{1}{n})$, one has $nx(1-x)^{n-1}\leq n\cdot\frac{1}{n}\cdot(1-\frac{1}{n})^{n-1}=(1-\frac{1}{n})^{n-1}\leq1;$\\
For $x\in[\frac{1}{n},1-\frac{1}{n}],$ let $r(x)=x(1-x)^{n-1},$ we can obtain $r'(x)=0,$ then $x=\frac{1}{n},$ one has \\
$nx(1-x)^{n-1}\leq n\cdot \frac{1}{n}\cdot(1-\frac{1}{n})^{n-1}=(1-\frac{1}{n})^{n-1}\leq1;$\\
For $x\in(1-\frac{1}{n},1], $  one has $nx(1-x)^{n-1}\leq n\cdot(1-x)^{n-1}\leq n\cdot(\frac{1}{n})^{n-1}=\frac{1}{n^{n-2}}\leq1.$\\
Hence
$$\varphi^{\lambda}(x)|H_{1}|\leq4\alpha\|\varphi^{\lambda} f'\|.\eqno{(11)}$$
By the similar method, we can get
$$\varphi^{\lambda}(x)|H_{3}|\leq4\alpha\|\varphi^{\lambda} f'\|.\eqno{(12)}$$
Next,
\begin{eqnarray*}
H_{2}&=&\sum^{n-1}_{k=1}\frac{\int^{1} _{0}t^{(n-1)k-1}(1-t)^{(n-1)(n-1-k)-1}}{B((n-1)k,(n-1)(n-1-k))}\left[f(\beta t+(1-\beta)\frac{k}{n-1})-f(x)\right]dt\\
&&\times\left[J_{n,k}^{\alpha}(x)-J_{n,k+1}^{\alpha}(x)\right]'\\
&=&\sum^{n-1}_{k=1}\frac{\int^{1} _{0}t^{(n-1)k-1}(1-t)^{(n-1)(n-1-k)-1}}{B((n-1)k,(n-1)(n-1-k))}\left[\int^{\beta t+(1-\beta)\frac{k}{n-1}}_{x}f'(u)du\right]dt\left[J_{n,k}^{\alpha}(x)-J_{n,k+1}^{\alpha}(x)\right]',
\end{eqnarray*}
from (2), we have
\begin{eqnarray*}
\varphi^{\lambda}(x)\cdot\left|\int^{\beta t+(1-\beta)\frac{k}{n-1}}_{x}f'(u)du\right|
&\leq& \varphi^{\lambda}(x)\cdot\|\varphi^{\lambda} f'\|\cdot \left|\int^{\beta t+(1-\beta)\frac{k}{n-1}}_{x}\frac{1}{\varphi^{\lambda}(u)}du\right|\\
&\leq& 4\|\varphi^{\lambda} f'\|\cdot \left|\beta t+(1-\beta)\frac{k}{n-1}-x\right|,
\end{eqnarray*}
then,
\begin{eqnarray*}
&&\varphi^{\lambda}(x)|H_{2}|\\
&\leq&4\|\varphi^{\lambda} f'\|\sum^{n-1}_{k=1}\frac{\int^{1} _{0}t^{(n-1)k-1}(1-t)^{(n-1)(n-1-k)-1}}{B((n-1)k,(n-1)(n-1-k))}\left|\beta t+(1-\beta)\frac{k}{n-1}-x\right|dt\\
&&\times\left|\left[J_{n,k}^{\alpha}(x)-J_{n,k+1}^{\alpha}(x)\right]'\right|\\
&=&4\alpha\|\varphi^{\lambda} f'\|\sum^{n-1}_{k=1}\frac{\int^{1} _{0}t^{(n-1)k-1}(1-t)^{(n-1)(n-1-k)-1}}{B((n-1)k,(n-1)(n-1-k))}\left|\beta t+(1-\beta)\frac{k}{n-1}-x\right|dt\\
&&\times\left|J_{n,k}^{\alpha-1}(x)\cdot[J_{n,k+1}'(x)+p_{n,k}'(x)]-J_{n,k+1}^{\alpha-1}(x)\cdot J_{n,k+1}'(x)\right|\\
&=&4\alpha\|\varphi^{\lambda} f'\|\left\{\sum^{n-1}_{k=1}\frac{\int^{1} _{0}t^{(n-1)k-1}(1-t)^{(n-1)(n-1-k)-1}}{B((n-1)k,(n-1)(n-1-k))}\left|\beta t+(1-\beta)\frac{k}{n-1}-x\right|dt\right.\\
&&\times\left|J_{n,k}^{\alpha-1}(x)-J_{n,k+1}^{\alpha-1}(x)\right|J_{n,k+1}'(x)\\
&&+\left.\sum^{n-1}_{k=1}\frac{\int^{1} _{0}t^{(n-1)k-1}(1-t)^{(n-1)(n-1-k)-1}}{B((n-1)k,(n-1)(n-1-k))}\left|\beta t+(1-\beta)\frac{k}{n-1}-x\right|dtJ_{n,k}^{\alpha-1}(x)|p_{n,k}'(x)|\right\},
\end{eqnarray*}
write $$
\varphi^{\lambda}(x)|H_{2}|\leq4\alpha\|\varphi^{\lambda} f'\|\cdot(A+B),\eqno{(13)}$$
we will estimate A and B on $E_n^C$ and $E_n$ respectively.\\
(I). For $x\in E_{n}^{c}=[0,\frac{1}{n})\bigcup(1-\frac{1}{n},1]$:
\begin{eqnarray*}
B&\leq&\sum^{n-1}_{k=1}\frac{\int^{1} _{0}t^{(n-1)k-1}(1-t)^{(n-1)(n-1-k)-1}}{B((n-1)k,(n-1)(n-1-k))}\left|\beta t+(1-\beta)\frac{k}{n-1}-x\right|dt\\
&&\times n|p_{n-1,k-1}(x)-p_{n-1,k}(x)|\\
&\leq&\sum^{n-1}_{k=1}\frac{\int^{1} _{0}t^{(n-1)k-1}(1-t)^{(n-1)(n-1-k)-1}}{B((n-1)k,(n-1)(n-1-k))}\left|\beta t+(1-\beta)\frac{k}{n-1}-x\right|dt\cdot n \cdot p_{n-1,k-1}(x)\\
&&+\sum^{n-1}_{k=1}\frac{\int^{1} _{0}t^{(n-1)k-1}(1-t)^{(n-1)(n-1-k)-1}}{B((n-1)k,(n-1)(n-1-k))}\left|\beta t+(1-\beta)\frac{k}{n-1}-x\right|dt\cdot n\cdot p_{n-1,k}(x),
\end{eqnarray*}
we write $$
B\leq L_{1}+L_{2}.\eqno{(14)}$$
Since $E_{n-1,\beta}((t-x)^{2};x)\leq \frac{2\varphi^{2}(x)}{n-1},$ by the Cauchy-Schwarz inequality, we get
\begin{eqnarray*}
L_{2}&=&n\left(\sum^{n-1}_{k=1}F_{n-1,k}^{(\beta)}((t-x)^{2};x)\cdot p_{n-1,k}(x)\right)^{\frac{1}{2}}\cdot\left(\sum^{n-1}_{k=1}p_{n-1,k}(x)\right)^{\frac{1}{2}}\\
&\leq& n \left(E_{n-1,\beta}((t-x)^{2};x)\right)^{\frac{1}{2}}\cdot 1\\
&\leq& 2n\cdot \frac{\sqrt{x(1-x)}}{\sqrt{n}}.
\end{eqnarray*}
For $x\in [0,\frac{1}{n}),$
$$L_{2}\leq 2n\cdot \frac{\sqrt{x}\cdot\sqrt{1-x}}{\sqrt{n}}\leq 2n\cdot \frac{\sqrt{x}}{\sqrt{n}}\leq 2;\eqno{(15)}$$
for $x\in (1-\frac{1}{n},1],$
$$L_{2}\leq 2n\cdot \frac{\sqrt{x}\cdot\sqrt{1-x}}{\sqrt{n}}\leq 2n\cdot \frac{\sqrt{1-x}}{\sqrt{n}}\leq 2.\eqno{(16)}$$
Next, we will estimate $L_{1}$:
\begin{eqnarray*}
L_{1}
&\leq& \sum^{n-2}_{j=0}\frac{\int^{1} _{0}t^{(n-1)(j+1)-1}(1-t)^{(n-1)[n-1-(j+1)]-1}}{B((n-1)(j+1),(n-1)[n-1-(j+1)])}\left|\beta t+(1-\beta)\frac{j+1}{n-1}-x\right|dt\cdot np_{n-1,j}(x)\\
&\leq& \frac{\int^{1} _{0}t^{(n-1)-1}(1-t)^{(n-1)(n-1-1)-1}}{B((n-1),(n-1)(n-1-1))}\left|\beta t+(1-\beta)\frac{1}{n-1}-x\right|dt\cdot n\cdot p_{n-1,0}(x)\\
&&+ \sum^{n-2}_{j=1}\frac{\int^{1} _{0}t^{(n-1)(j+1)-1}(1-t)^{(n-1)[n-1-(j+1)]-1}}{B((n-1)(j+1),(n-1)[n-1-(j+1)])}\left|\beta t+(1-\beta)\frac{j}{n-1}-x\right|dt\cdot np_{n-1,j}(x)\\
&&+\sum^{n-2}_{j=1}\frac{\int^{1} _{0}t^{(n-1)(j+1)-1}(1-t)^{(n-1)[n-1-(j+1)]-1}}{B((n-1)(j+1),(n-1)[n-1-(j+1)])}\cdot\frac{1-\beta}{n-1}dt\cdot n\cdot p_{n-1,j}(x)\\
&\leq&F_{n-1,1}^{(\beta)}(t;x)\cdot n\cdot p_{n-1,0}(x)+x\cdot n\cdot p_{n-1,0}(x)+L_{3}+\sum^{n-2}_{j=1}\frac{1-\beta}{n-1}dt\cdot n\cdot p_{n-1,j}(x)\\
&\leq& \frac{1}{n-1}\cdot n\cdot (1-x)^{n-1}+x\cdot n\cdot (1-x)^{n-1}+L_{3}+\frac{n}{n-1}\cdot (1-\beta)\sum^{n-2}_{j=1}p_{n-1,j}(x)\\
&\leq& 4+n\cdot x(1-x)^{n-1}+L_{3},
\end{eqnarray*}
For $x\in [0,\frac{1}{n}),$ $n\cdot x(1-x)^{n-1}\leq n\cdot x\leq1;$ for $x\in (1-\frac{1}{n},1],$ $n\cdot x(1-x)^{n-1}\leq n\cdot (1-x)^{n-1}\leq \frac{1}{n^{n-2}}\leq 1$,
using the fact, $B(p,q)=\frac{p-1}{q}B(p-1,q+1) \hspace{0.5cm}(p,q>0)$, and for $j\geq1$
\begin{eqnarray*}
&&\int^{1} _{0}t^{(n-1)j-1+(n-1)}(1-t)^{(n-1)(n-1-j)-1-(n-1)}\cdot|\beta t+(1-\beta)\frac{j}{n-1}-x|dt\\
&\leq&\frac{(n-1)j+(n-1)}{(n-1)(n-1-j)-(n-1)}\cdot \frac{(n-1)j+(n-2)}{(n-1)(n-1-j)-(n-2)}\cdot\cdot\cdot\frac{(n-1)j+1}{(n-1)(n-1-j)-1}\\
&&\times[\int^{1} _{0}t^{(n-1)j-1}(1-t)^{(n-1)(n-1-j)-1}\cdot|\beta t+(1-\beta)\frac{j}{n-1}-x|dt\\
&&+\frac{n}{(n-1)^{2}}\cdot \beta\cdot\int^{1} _{0}t^{(n-1)j-1}(1-t)^{(n-1)(n-1-j)-1}dt],
\end{eqnarray*}
we can obtain
\begin{eqnarray*}
L_{3}
&\leq& \sum^{n-2}_{j=1}\frac{\int^{1} _{0}t^{(n-1)j-1}(1-t)^{(n-1)(n-1-j)-1}}{B((n-1)j,(n-1)(n-1-j))}\left|\beta t+(1-\beta)\frac{j}{n-1}-x\right|dt\cdot n\cdot p_{n-1,j}(x)\\
&&+\sum^{n-2}_{j=1}\frac{\int^{1} _{0}t^{(n-1)j-1}(1-t)^{(n-1)(n-1-j)-1}}{B((n-1)j,(n-1)(n-1-j))}dt\cdot\frac{n}{(n-1)^{2}}\cdot \beta\cdot p_{n-1,j}(x)\\
&\leq& L_{2}+\sum^{n-2}_{j=1}\frac{\int^{1} _{0}t^{(n-1)j-1}(1-t)^{(n-1)(n-1-j)-1}}{B((n-1)j,(n-1)(n-1-j))}dt\cdot \beta\cdot p_{n-1,j}(x)\\
&\leq& L_{2}+\beta\cdot\sum^{n-2}_{j=1} p_{n-1,j}(x)\\
&\leq& L_{2}+\beta\leq 3,
\end{eqnarray*}
so  $$L_{1}\leq L_{3}+5\leq 8,\eqno{(17)}$$
from (14)-(17), we know $B\leq 10.$\\
Noting that $J_{n,0}'(x)=0,$ and  $x\in E_n^C$, we get
\begin{eqnarray*}
A&=&\sum^{n-1}_{k=1}\frac{\int^{1} _{0}t^{(n-1)k-1}(1-t)^{(n-1)(n-1-k)-1}}{B((n-1)k,(n-1)(n-1-k))}\left|\beta t+(1-\beta)\frac{k}{n-1}-x\right|dt\cdot J_{n,k}^{\alpha-1}(x)J_{n,k}'(x)\\
&&-\sum^{n-1}_{k=1}\frac{\int^{1} _{0}t^{(n-1)k-1}(1-t)^{(n-1)(n-1-k)-1}}{B((n-1)k,(n-1)(n-1-k))}\left|\beta t+(1-\beta)\frac{k}{n-1}-x\right|dt\cdot J_{n,k}^{\alpha-1}(x)p_{n,k}'(x)\\
&&-\sum^{n-1}_{k=1}\frac{\int^{1} _{0}t^{(n-1)k-1}(1-t)^{(n-1)(n-1-k)-1}}{B((n-1)k,(n-1)(n-1-k))}\left|\beta t+(1-\beta)\frac{k}{n-1}-x\right|dt\cdot J_{n,k+1}^{\alpha-1}(x)J_{n,k+1}'(x)\\
&\leq&L_{4}+B-\sum^{n-1}_{k=1}\frac{\int^{1} _{0}t^{(n-1)k-1}(1-t)^{(n-1)(n-1-k)-1}}{B((n-1)k,(n-1)(n-1-k))}\left|\beta t+(1-\beta)\frac{k}{n-1}-x\right|dt\\
&&\times J_{n,k+1}^{\alpha-1}(x)J_{n,k+1}'(x),
\end{eqnarray*}
let $k=j+1$, 
\begin{eqnarray*}
L_{4}&=&\frac{\int^{1} _{0}t^{(n-1)-1}(1-t)^{(n-1)(n-1-1)-1}}{B((n-1),(n-1)(n-1-1))}\left|\beta t+(1-\beta)\frac{1}{n-1}-x\right|dtJ_{n,1}^{\alpha-1}(x)J_{n,1}'(x)\\
&&+\sum^{n-2}_{j=1}\frac{\int^{1} _{0}t^{(n-1)(j+1)-1}(1-t)^{(n-1)[n-1-(j+1)]-1}}{B((n-1)(j+1),(n-1)[n-1-(j+1)])}\left|\beta t+(1-\beta)\frac{j+1}{n-1}-x\right|dt\\
&&\times J_{n,j+1}^{\alpha-1}(x)J_{n,j+1}'(x)\\
&\leq&[F_{n-1,1}^{(\beta)}(t;x)+x]J_{n,1}^{\alpha-1}(x)J_{n,1}'(x)\\
&&+\sum^{n-2}_{j=1}\frac{\int^{1} _{0}t^{(n-1)(j+1)-1}(1-t)^{(n-1)[n-1-(j+1)]-1}}{B((n-1)(j+1),(n-1)[n-1-(j+1)])}\left|\beta t+(1-\beta)\frac{j}{n-1}-x\right|dt\\
&&\times J_{n,j+1}^{\alpha-1}(x)J_{n,j+1}'(x)\\
&&+\sum^{n-2}_{j=1}\frac{\int^{1} _{0}t^{(n-1)(j+1)-1}(1-t)^{(n-1)[n-1-(j+1)]-1}}{B((n-1)(j+1),(n-1)[n-1-(j+1)])}\left|\frac{1-\beta}{n-1}\right|dt\cdot J_{n,j+1}^{\alpha-1}(x)J_{n,j+1}'(x)\\
&\leq& \frac{n}{n-1}\cdot (1-x)^{n-1}+nx\cdot (1-x)^{n-1}+\frac{n-2}{n-1}\\
&&+\sum^{n-2}_{j=1}\frac{\int^{1} _{0}t^{(n-1)(j+1)-1}(1-t)^{(n-1)[n-1-(j+1)]-1}}{B((n-1)(j+1),(n-1)[n-1-(j+1)])}\left|\beta t+(1-\beta)\frac{j}{n-1}-x\right|dt\\
&&\times J_{n,j+1}^{\alpha-1}(x)J_{n,j+1}'(x)\\
&\leq& 4+ \sum^{n-2}_{j=1}\frac{\int^{1} _{0}t^{(n-1)j-1}(1-t)^{(n-1)(n-1-j)-1}}{B((n-1)j,(n-1)(n-1-j))}\left|\beta t+(1-\beta)\frac{j}{n-1}-x\right|dtJ_{n,j+1}^{\alpha-1}(x)J_{n,j+1}'(x),
\end{eqnarray*}
hence $A\leq 10+B\leq 14. $\\
Combing A and B, we have,
$$|\varphi^{\lambda}(x)H_{2}|\leq 96\alpha\|\varphi^{\lambda} f'\|.\eqno{(18)}$$
From (11), (12) and (18), for $x\in E_{n}^{C}$, we have
$$|\varphi^{\lambda}(x)\cdot(L_{n,\beta}^{(\alpha)}(f;x))'|\leq 104\alpha\|\varphi^{\lambda}f'\| .$$
(II). $x\in E_{n}=[\frac{1}{n}, 1-\frac{1}{n}]$: Noting $p_{n,k}'(x)=\frac{n}{\varphi^{2}(x)}(\frac{k}{n}-x)\cdot p_{n,k}(x)$, using the  Cauchy-Schwarz inequality,
\begin{eqnarray*}
B&\leq& \sum^{n-1}_{k=1}\frac{\int^{1} _{0}t^{(n-1)k-1}(1-t)^{(n-1)(n-1-k)-1}}{B((n-1)k,(n-1)(n-1-k))}\left|\beta t+(1-\beta)\frac{k}{n-1}-x\right|dt\cdot |p_{n,k}'(x)|\\
&\leq& \left(\sum^{n-1}_{k=1}\frac{\int^{1} _{0}t^{(n-1)k-1}(1-t)^{(n-1)(n-1-k)-1}}{B((n-1)k,(n-1)(n-1-k))}\left(\beta t+(1-\beta)\frac{k}{n-1}-x\right)^{2}dt\cdot p_{n,k}(x)\right)^{\frac{1}{2}}\\
&&\times\left(\sum^{n-1}_{k=1}\frac{\int^{1} _{0}t^{(n-1)k-1}(1-t)^{(n-1)(n-1-k)-1}dt}{B((n-1)k,(n-1)(n-1-k))}\cdot p_{n,k}(x)\cdot(\frac{k}{n}-x)^{2}\right)^{\frac{1}{2}}\cdot \frac{n}{\varphi^{2}(x)}\\
&\leq&\left(\sum^{n-1}_{k=1}p_{n,k}(x)F_{n-1,k}^{(\beta)}((t-x)^{2};x)\right)^{\frac{1}{2}} \cdot\left(B_{n}((t-x)^{2};x)\right)^{\frac{1}{2}}\cdot\frac{n}{\varphi^{2}(x)}\\
&\leq&\left(L_{n,\beta}((t-x)^{2};x)\right)^{\frac{1}{2}}\cdot \left(B_{n}((t-x)^{2};x)\right)^{\frac{1}{2}}\cdot\frac{n}{\varphi^{2}(x)}\\
&\leq& 3\cdot\frac{\varphi(x)}{\sqrt{n-1}}\cdot \frac{\varphi(x)}{\sqrt{n}}\cdot \frac{n}{\varphi^{2}(x)}\leq 6\  \    (for\  \
   n\geq 2),
\end{eqnarray*}
in the last inequality, we have used Lemma 7(2).\\
Using the same method as $x\in  E_{n}^{C},$ we get $A\leq 4+B\leq 10,$ then
$$|\varphi^{\lambda}(x)H_{2}|\leq 64\alpha\|\varphi ^{\lambda}f'\|.$$
Hence, for $x\in E_{n}$, we have
$$|\varphi^{\lambda}(x)\cdot(L_{n,\beta}^{(\alpha)}(f;x))'|\leq|\varphi^{\lambda}(x)H_{1}|+|\varphi^{\lambda}(x)H_{2}|+|\varphi^{\lambda}(x)H_{3}| \leq 72\alpha\|\varphi ^{\lambda}f'\| .$$
Then, we can get the desired result.\\
{\bf Acknowledgments}\\
This work was supported by NSF of China [Grant No. 11871191, 11571089],
 NSF of Hebei Province [Grant No. ZD2019053] and postgraduate Innovation Funding Project of Hebei Normal University(Grant No. CXZZSS2019059).\\
We express our gratitude to the referees for their helpful suggestions.

\end{document}